# A GENERALIZED ALLWRIGHT FORMULA AND THE VECTOR RICCATI EQUATION

KURT MUNK ANDERSEN AND ALLAN SANDQVIST

### ABSTRACT

Let $\varphi(t,\tau,\xi)$ denote the general solution of a scalar differential equation $\frac{dx}{dt} = f(t,x)$. A formula of Allwright involving the derivative $\frac{\partial^3}{\partial \xi^3}\varphi(t,\tau,\xi)$ is generalized to the case of a differential system by means of the Kronecker product. The Allwright formula is connected with the Riccati equation $\frac{dx}{dt} = a(t) + b(t)x + c(t)x^2$, and in a similar way the generalized formula is connected with the vector Riccati equation $\frac{d\underline{x}}{dt} = \underline{a}(t) + \underline{\underline{B}}(t)\underline{x} + \left(\underline{c}^T(t)\underline{x}\right)\underline{x}$.

Moreover, the classical result that a scalar differential equation $\frac{dx}{dt} = f(t,x)$ is a Riccati equation if and only if the general solution $\varphi(t,\tau,\xi)$ is a fractional linear function in $\xi$ is generalized to the case of a differential system.

## 1. Introduction

Consider a scalar differential equation

$$\frac{dx}{dt} = f(t,x)\,, (t,x) \in I \times R \tag{1}$$

where $I \subseteq R$ is an open interval and $f(t,x), (t,x) \in I \times R$ is a continuous function, which is $C^n$ in $x$ (for a suitable n). Let $\varphi(t,\tau,\xi), t \in I(\tau,\xi)$ denote the maximal solution of (1) through the point $(\tau,\xi) \in I \times R$. The following formulae are well known (a prime indicates differention with respect to $\xi$ or $x$):

$$\varphi'(t,\tau,\xi) = e^{\int_\tau^t f'(s,\varphi(s,\tau,\xi))ds} \tag{2}$$

---

$$\varphi''(t,\tau,\xi) = e^{\int_\tau^t f'(s,\varphi(s,\tau,\xi))ds} \int_\tau^t f''(s,\varphi(s,\tau,\xi)) \, e^{\int_\tau^s f'(u,\varphi(u,\tau,\xi))du} ds \qquad (3)$$

$$2\varphi'(t,\tau,\xi)\varphi'''(t,\tau,\xi) - 3(\varphi''(t,\tau,\xi))^2 = 2\int_\tau^t f'''(s,\tau,\xi) \, e^{2\int_\tau^s f'(u,\varphi(u,\tau,\xi))du} ds \qquad (4)$$

The two first formulae are classical. The third formula is shown by Allwright (see [1, Theorem 2]). The following question arises: is it possible to generalize these formulae to a system of n equations

$$\frac{d\underline{x}}{dt} = \underline{f}(t,\underline{x}), (t,\underline{x}) \in I \times R^n \qquad (5)$$

where $\underline{f}(t,\underline{x}), (t,\underline{x}) \in I \times R^n$ is a continuous function, which is $C^n$ in $\underline{x}$ (for a suitable n)? It is classical that the answer is affirmative as to (2) (see, e.g., [5]). In [**7**] it is proved that it is also affirmative as to (3). In section 4 in the present paper we prove that the Allwright formula (4) can be generalized, too. In the following $\underline{D}$ indicates the Jacobi matrix with respect to $\underline{x}$ (or $\underline{\xi}$). Let $\underline{\varphi}(t,\tau,\underline{\xi}), t \in I(\tau,\underline{\xi})$ denote the maximal solution of (5) through the point $(\tau,\underline{\xi}) \in I \times R^n$ and let $\underline{\Phi}(t,\tau,\underline{\xi}), t \in I(\tau,\underline{\xi})$ denote the fundamental matrix of the first variation

$$\frac{d\underline{z}}{dt} = \underline{D}\underline{f}(t,\underline{\varphi}(t,\tau,\underline{\xi}))\underline{z}, t \in I(\tau,\underline{\xi}) \qquad (6)$$

for which $\underline{\Phi}(\tau,\tau,\underline{\xi}) = \underline{E}$, the $n \times n$ unit matrix. The generalized version of (2) then is

$$\underline{D}\underline{\varphi}(t,\tau,\underline{\xi}) = \underline{\Phi}(t,\tau,\underline{\xi}) \qquad (7)$$

The generalized version of (3) is

$$\underline{D}^2\underline{\varphi}(t,\tau,\underline{\xi}) = \underline{\Phi}(t,\tau,\underline{\xi}) \int_\tau^t \underline{\Phi}^{-1}(s,\tau,\underline{\xi}) \underline{D}^2\underline{f}(t,\underline{\varphi}(s,\tau,\underline{\xi}))\underline{\Phi}(s,\tau,\underline{\xi}) \otimes \underline{\Phi}(s,\tau,\underline{\xi}) \, ds \qquad (8)$$

Here $\underline{D}^2\underline{\varphi}$ denotes the block matrix $\left[\frac{\partial}{\partial\xi_1}\underline{D}\underline{\varphi}, ..., \frac{\partial}{\partial\xi_n}\underline{D}\underline{\varphi}\right]$, and analogous as to $\underline{D}^2\underline{f}$, and $\otimes$ denotes the Kronecker product. In view of (7) the matrix $\underline{\Phi}$ can be replaced by the matrix $\underline{D}\underline{\varphi}$ in (8). We will prove the following formula, where the arguments are omitted for the sake of clearness:

$$\underline{d}^3\,\underline{\varphi}\otimes\underline{d\varphi}+\underline{d\varphi}\otimes\underline{d}^3\,\underline{\varphi}-3\big[\underline{d}^2\underline{\varphi}\big]^2=\underline{\underline{D}\varphi}\otimes\underline{\underline{D}\varphi}\bigg\{\int_\tau^t\big[\underline{\underline{D}\varphi}\otimes\underline{\underline{D}\varphi}\big]^{-1}\big(\underline{\underline{D}}^3\,\underline{f}\otimes\underline{\underline{E}}+\underline{\underline{E}}\otimes\underline{\underline{D}}^3\,\underline{f}\big)\big(\underline{d\varphi}\big)^4\,ds\bigg\}$$

$$+3\bigg\{\int_\tau^t\big[\underline{\underline{D}\varphi}\otimes\underline{\underline{D}\varphi}\big]^{-1}\big(\underline{\underline{D}}^2\,\underline{f}\otimes\underline{\underline{E}}-\underline{\underline{E}}\otimes\underline{\underline{D}}^2\,\underline{f}\big)\big(\underline{d}^2\,\underline{\varphi}\otimes\big(\underline{d\varphi}\big)^2-\big(\underline{d\varphi}\big)^2\otimes\underline{d}^2\,\underline{\varphi}\big)ds\bigg\}\qquad(9)$$

In this formula $\underline{d}^k\underline{\varphi}$ denotes the k'th differential of $\underline{\varphi}$ with respect to $\underline{\xi}$ and

$\underline{\underline{D}}^3\,\underline{f}$ denotes the block matrix $\left[\dfrac{\partial}{\partial x_1}\underline{\underline{D}}^2\,\underline{f},...,\dfrac{\partial}{\partial x_n}\underline{\underline{D}}^2\,\underline{f}\right]$. Moreover, $\big(\underline{d\varphi}\big)^k$ denotes the Kronecker produkt of $\underline{d\varphi}$ k times. In the integrand $\underline{\underline{D}}^k\,\underline{f}$ shall be taken in the point $\big(s,\underline{\varphi}(s,\tau,\underline{\xi})\big)$. The reason why differentials occur will be

motivated in connection with the proof. If n=1 then (9) becomes the

Allwright formula (4), because the second integral vanishes. If n≥2 this integral normally does not vanish, because the factors in a Kronecker product in general do not commute, so that the two last factors in the integrand are not necesssarily zero. The Allwright formula (4) is proved by differention of (3).

In principle formula (9) is proved similarly. However, as is to be expected from the complexity of formula (9) this requires a great deal of matrix technique and auxiliary propositions.

In section 5 we generalize a well known result connected with the Allwright

formula. As shown in [2] this formula can be written

$$\{\varphi,\xi\}=\int_\tau^t f'''(s,\varphi(s,\tau,\xi))e^{2\int_\tau^s f'(u,\varphi(u,\tau,\varphi))du}\,ds\qquad(10)$$

Here $\{\varphi,\xi\}$ denotes the Schwarzian derivative of $\varphi$ with respect to $\xi$. With

omitted arguments this is:

$$\{\varphi,\xi\}=\frac{\varphi'''}{\varphi'}-\frac{3}{2}\left\{\frac{\varphi''}{\varphi'}\right\}^2\qquad(11)$$

It is classical that $\{\varphi,\xi\}$ is identically zero ( i.e. the left hand side of (4) is

identically zero) if and only if (1) is a Riccati equation

$$\frac{dx}{dt}=a(t)+b(t)x+c(t)x^2\ ,\ t\in I\qquad(12)$$

see [4, pp. 274,647-649,671]. This result is also a direct consequence of (10).

It turns out that similarly the expression on the left hand side of (9) is

identically zero if and only if system (5) has the form

$$\frac{d\underline{x}}{dt} = \underline{a}(t) + \underline{\underline{B}}(t)\underline{x} + \left(\underline{c}^T(t)\underline{x}\right)\underline{x}, \, t \in I \tag{13}$$

where $\underline{c}^T$ denotes the transpose of $\underline{c}$. If n=1 equation (13) becomes a Riccati equation. Hence, we call (13) a vector Riccati equation. In section 6 we deal with another well known result, namely that (1) is a Riccati equation if and only if the general solution has the form

$$\varphi(t, \tau, \xi) = \frac{\alpha(t, \tau)\xi + \beta(t, \tau)}{\gamma(t, \tau)\xi + \delta(t, \tau)} \tag{14}$$

i.e., a fractional linear function in $\xi$ (see [4] loc. cit.). We prove that (5) is a vector Riccati equation if and only if the general solution has the form

$$\underline{\varphi}(t, \tau, \underline{\xi}) = \frac{\underline{\underline{A}}(t, \tau)\underline{\xi} + \underline{\beta}(t, \tau)}{\underline{\gamma}^T(t, \tau)\underline{\xi} + \delta(t, \tau)} \tag{15}$$

which we call a fractional linear vector function in $\underline{\xi}$. If n=1 this becomes the expression (14). Sections 2 and 3 are auxiliary in character.

## 2. Some matrix notions

In the following vectors are writtes as columns. The usual basis for $R^n$ is denoted $\underline{e}_1, \underline{e}_2, ..., \underline{e}_n$. The n×n unit matrix is denoted $\underline{\underline{E}}_n$, where the subscript is omitted if the dimension n is obvious.

The Kronecker product of an n×m matrix $\underline{\underline{A}}$ and a p×q matrix $\underline{\underline{B}}$ is defined as the np×mq matrix

$$\underline{\underline{A}} \otimes \underline{\underline{B}} = [a_{rs}\underline{\underline{B}}], \, r = 1, 2, ..., p, \, s = 1, 2, ..., q$$

The usual matrix product and the Kronecker product are connected via the *product rule*

$$\left(\underline{\underline{A}} \otimes \underline{\underline{B}}\right)\left(\underline{\underline{C}} \otimes \underline{\underline{D}}\right) = \left(\underline{\underline{A}}\underline{\underline{C}}\right) \otimes \left(\underline{\underline{B}}\underline{\underline{D}}\right) \tag{16}$$

which holds provided the dimensions of the matrices are such that the various expressions exist. This rule and other properties of the Kronecker product can be found in [3], see in particular the table on page 122. The Kronecker product of a matrix $\underline{\underline{A}}$ with itself p times is denoted $\underline{\underline{A}}^{\otimes p}$. If $\underline{\underline{A}}$ is a column $\underline{h}$ we omit the $\otimes$ before p and write $\underline{h}^p$ (the scalar product of $\underline{h}$ with itself will not occur, so there will be no confusion for p=2).

Let $\Omega \subseteq R^n$ be an open domain. If $\underline{\underline{M}}(\underline{x})$, $\underline{x} \in \Omega$ is a $C^1$ matrix function we write

$$\underline{\underline{DM}} = \left[ \frac{\partial}{\partial x_1} \underline{\underline{M}} \; \frac{\partial}{\partial x_2} \underline{\underline{M}} \; ... \; \frac{\partial}{\partial x_n} \underline{\underline{M}} \right] \qquad (17)$$

If $\underline{\underline{M}}$ is $1 \times 1$ then (17) is the usual gradient. Let $\underline{g} : \Omega \to R^p$ be a $C^q$ vector function. Then (17) becomes the functional matrix $\underline{\underline{Dg}}$, if $\underline{g}$ is inserted for $\underline{\underline{M}}$. We write $\underline{\underline{D}}(\underline{\underline{Dg}}) = \underline{\underline{D}}^2 \underline{g}$, and recursively $\underline{\underline{D}}(\underline{\underline{D}}^{m-1} \underline{g}) = \underline{\underline{D}}^m \underline{g}$, $m = 2,...,q$. The first differential of $\underline{g}$ is

$$\underline{dg}(\underline{x},\underline{h}) = \sum_{k=1}^{n} \frac{\partial \underline{g}(\underline{x})}{\partial x_k} h_k = \underline{\underline{Dg}}(\underline{x})\underline{h} \; , \; \underline{x} \in \Omega \, , \; \underline{h} \in R^n$$

The second differential is

$$\underline{d}^2 \underline{g}(\underline{x},\underline{h}) = \sum_{k,l=1}^{n} \frac{\partial^2 \underline{g}(\underline{x})}{\partial x_k \partial x_l} h_k h_l = \sum_{l=1}^{n} \frac{\partial}{\partial x_l} (\underline{\underline{Dg}}(\underline{x})\underline{h}) h_l = \sum_{l=1}^{n} \frac{\partial}{\partial x_l} \underline{\underline{Dg}}(\underline{x})(h_l \underline{h})$$

$$= \left[ \frac{\partial}{\partial x_1} \underline{\underline{Dg}}(\underline{x}) \; ... \; \frac{\partial}{\partial x_n} \underline{\underline{Dg}}(\underline{x}) \right] \underline{h} \otimes \underline{h} = \underline{\underline{D}}^2 \underline{g}(\underline{x})\underline{h}^2 \, , \; \underline{x} \in \Omega \, , \; \underline{h} \in R^n$$

and in the same way the m'th differential can be written

$$\underline{d}^m \underline{g}(\underline{x},\underline{h}) = \sum_{k_1,k_2,...,k_n=1}^{n} \frac{\partial^m \underline{g}(\underline{x})}{\partial x_{k_1} \partial x_{k_2} ... \partial x_{k_n}} h_{k_1} h_{k_2} ... h_{k_n} = \underline{\underline{D}}^m \underline{g}(\underline{x})\underline{h}^m \, , \; \underline{x} \in \Omega \, , \; \underline{h} \in R^n \, ,$$

m=1, 2, ... , q. Hence Taylor's formula can be written

$$\underline{g}(\underline{x} + \underline{h}) - \underline{g}(\underline{x}) = \sum_{k=1}^{q-1} \frac{1}{k!} \underline{\underline{D}}^k \underline{g}(\underline{x})\underline{h}^k + \frac{1}{q!} \underline{\underline{D}}^q \underline{g}(\underline{x} + \Theta \underline{h})\underline{h}^n \qquad (18)$$

According to (17) the matrix $\underline{\underline{D}}^2 \underline{g}$ consists of n n×n blocks

$\frac{\partial}{\partial x_k} \underline{\underline{Dg}}, k = 1,2,...,n$. The i'th column in $\frac{\partial}{\partial x_j} \underline{\underline{Dg}}$ and the j'th column

in $\frac{\partial}{\partial x_i} \underline{\underline{Dg}}$ are equal, namely the second order derivative $\frac{\partial^2 \underline{g}}{\partial x_i \partial x_j}$. An n×n² block matrix $\underline{\underline{A}} = \left[ \underline{\underline{A}}_1 \; \underline{\underline{A}}_2 \; ... \; \underline{\underline{A}}_n \right]$ with this property will be called *column symmetric* (briefly *c-symmetric*). So, the requirement is that the i'th column in $\underline{\underline{A}}_j$ equals the j'th column in $\underline{\underline{A}}_i$. This can be expressed

$$\underline{\underline{A}}(\underline{e}_i \otimes \underline{e}_j) = \underline{\underline{A}}(\underline{e}_j \otimes \underline{e}_i) \, , i,j = 1,2,...,n \qquad (19)$$

Observe namely that $\underline{e}_i \otimes \underline{e}_j = \begin{bmatrix} \underline{0} & \underline{0} & \dots & \underline{e}_j & \dots & \underline{0} \end{bmatrix}^T$ where $\underline{e}_j$ is at position i.

EXAMPLE 2.1. Let $\underline{a}$ be an n×1 column. Then the n×n² matrix

$\underline{a}^T \otimes \underline{\underline{E}} + \underline{\underline{E}} \otimes \underline{a}^T$ is c-symmetric. Indeed, by the product rule (16)

$$(\underline{a}^T \otimes \underline{\underline{E}} + \underline{\underline{E}} \otimes \underline{a}^T)(\underline{e}_i \otimes \underline{e}_j) = (\underline{a}^T \underline{e}_i) \otimes \underline{e}_j + \underline{e}_i \otimes (\underline{a}^T \underline{e}_j) = a_i \underline{e}_j + a_j \underline{e}_i$$

which remains unchanged if i and j are interchanged.

The definition easily can be extended. An m×n^q block matrix $\underline{\underline{M}}$ is called *c-symmetric* if $\underline{\underline{M}}\left(\underline{e}_{i_1} \otimes \underline{e}_{i_2} \otimes \dots \otimes \underline{e}_{i_q}\right)$ is the same column for all q! permutations of any fixed index set $(i_1, i_2, \dots, i_q)$. Consequently, all matrices $\underline{D}^m \underline{g}$ are c-symmetric. Note that q=1 is allowed, meaning that (trivially) any m×n matrix is c-symmetric. These matrices have the following important property.

THEOREM 2.1. *Let* $\underline{\underline{M}}$ *be an m×n^q c-symmetric matrix. Then*

$$\underline{\underline{M}} = \underline{0} \Leftrightarrow \forall \underline{x} \in R^n : \underline{\underline{M}} \underline{x}^q = \underline{0} \qquad (20)$$

*Proof.* For any $\underline{x} \in R^k$ we may write

$$\underline{x}^q = \left(\sum_{i=1}^n x_i \underline{e}_i\right)^q = \sum x_{i_1} x_{i_2} \dots x_{i_q} (\underline{e}_{i_1} \otimes \underline{e}_{i_2} \otimes \dots \otimes \underline{e}_{i_q})$$

where the last summation is over all index sets $(i_1, i_2, \dots, i_q)$. We get

$$\underline{\underline{M}} \underline{x}^q = q! \sum_{i_1 \le i_2 \le \dots i_q} x_{i_1} x_{i_2} \dots x_{i_q} \underline{\underline{M}}(\underline{e}_{i_1} \otimes \underline{e}_{i_2} \otimes \dots \otimes \underline{e}_{i_q})$$

since both $x_{i_1} x_{i_2} \dots x_{i_q}$ and $\underline{\underline{M}}(\underline{e}_{i_1} \otimes \underline{e}_{i_2} \otimes \dots \otimes \underline{e}_{i_q})$ are unchanged if the indices permute. From the identity theorem for polynomials in q variables it follows that $\underline{\underline{M}} \underline{x}^q = \underline{0}$ for all $\underline{x} \in R^n$ implies that $x_{i_1} x_{i_2} \dots x_{i_q} \underline{\underline{M}}(\underline{e}_{i_1} \otimes \underline{e}_2 \otimes \dots \otimes \underline{e}_{i_q}) = \underline{0}$ for all $i_1 \le i_2 \le \dots i_q$ and hence for all index sets $(i_1, i_2, \dots, i_q)$. Since the products $\underline{e}_{i_1} \otimes \underline{e}_{i_2} \otimes \dots \otimes \underline{e}_{i_q}$ form the usual basis for R^n we infer that $\underline{\underline{M}} = \underline{0}$. The opposite implication is trivial. This proves (20).

Note that $\left\{\underline{x}^q \,\middle|\, \underline{x} \in R^n\right\}$ is only a *proper* subset of R^n·, so the theorem is not trivial. The following two corollaries are immediate consequences of Theorem 2.1.

COROLLARY 2.1. *If* $\underline{\underline{A}}$ *is a nonsingular n×n matrix, then*

$$\underline{\underline{M}} = \underline{0} \Leftrightarrow \forall \underline{x} \in R^n : \underline{\underline{M}}(\underline{\underline{A}} \underline{x})^q = \underline{0} \qquad (21)$$

COROLLARY 2.2. *Let $\underline{M}_1$ and $\underline{M}_2$ be $m \times n^q$ c-symmetric matrices. Then*

$$\underline{M}_1 = \underline{M}_2 \Leftrightarrow \forall \underline{x} \in R^n : \underline{M}_1 \underline{x}^q = \underline{M}_2 \underline{x}^q \qquad (22)$$

It is obvious that a linear combination of $m \times n^q$ c-symmetric matrices is c-symmetric. If $\underline{A}$ is any $p \times m$ matrix and $\underline{M}$ is an $m \times n^q$ c-symmetric matrix, then $\underline{AM}$ is also c-symmetric.

## 3. Some propositions involving the Kronecker product

The proofs of the following propositions are placed in the appendix.

PROPOSITION 3.1. *Let $\underline{a}, \underline{b}$ and $\underline{c}$ be $n \times 1$ columns and let $\underline{C}$ be an $m \times n$ matrix. Then*

$$(\underline{a}^T \otimes \underline{b})\underline{c} = (\underline{a}^T \underline{c})\underline{b} \text{ and } (\underline{a} \otimes \underline{C})\underline{b} = \underline{a} \otimes (\underline{C}\underline{b}) \text{ and } (\underline{C} \otimes \underline{a})\underline{b} = (\underline{C}\underline{b}) \otimes \underline{a} \qquad (23)$$

PROPOSITION 3.2. *Let $\underline{M}$ be an $m \times n^2$ c-symmmetric matrix, $\underline{A}$ an $n \times p$ matrix and $\underline{b}$ an $n \times 1$ column. Then*

$$\underline{M}(\underline{A} \otimes \underline{b}) = \underline{M}(\underline{b} \otimes \underline{A}) \qquad (24)$$

COROLLARY 3.1. *Let $\underline{a}, \underline{b}$ and $\underline{c}$ be $n \times 1$ columns and $\underline{E} = \underline{E}_n$. Then*

$$(\underline{M} \otimes \underline{E})(\underline{a} \otimes \underline{b} \otimes \underline{c}) = (\underline{M} \otimes \underline{E})(\underline{b} \otimes \underline{a} \otimes \underline{c}) \qquad (25)$$

$$(\underline{E} \otimes \underline{M})(\underline{a} \otimes \underline{b} \otimes \underline{c}) = (\underline{E} \otimes \underline{M})(\underline{a} \otimes \underline{c} \otimes \underline{b}) \qquad (26)$$

*Proof.* Via the product rule we get $(\underline{M} \otimes \underline{E})(\underline{a} \otimes \underline{b} \otimes \underline{c}) = (\underline{M}(\underline{a} \otimes \underline{b})) \otimes \underline{c}$ and $(\underline{M} \otimes \underline{E})(\underline{b} \otimes \underline{a} \otimes \underline{c}) = (\underline{M}(\underline{b} \otimes \underline{a})) \otimes \underline{c}$. Then (25) follows from (24) with $\underline{A} = \underline{a}$, Rule (26) is proved analogously.

PROPOSITION 3.3. *Let $\underline{a}, \underline{b}$ and $\underline{c}$ be $n \times 1$ columns. Then*

$$\underline{a} \otimes \underline{b} + \underline{b} \otimes \underline{a} = \underline{0} \Rightarrow \underline{a} = \underline{0} \vee \underline{b} = \underline{0} \qquad (27)$$

$$\underline{a} \otimes \underline{b} + \underline{b} \otimes \underline{a} = \underline{c}^2 \Rightarrow \exists \alpha, \beta \in R : \underline{c} = \alpha \underline{a} = \beta \underline{b} \qquad (28)$$

*If $\underline{a} \neq \underline{0}$ then (28) can be sharpened :*

$$\underline{a} \otimes \underline{b} + \underline{b} \otimes \underline{a} = \underline{c}^2 \Rightarrow \exists \alpha \in R : \underline{c} = \alpha \underline{a} \wedge \underline{b} = \frac{1}{2}\alpha^2 \underline{a} \qquad (29)$$

For any matrices $\underline{A}$ and $\underline{B} = \left[ \underline{B}_1 \ \underline{B}_2 \ ... \underline{B}_n \right]$ we write for practical reasons

$$\underline{A} * \underline{B} = \left[ \underline{A} \otimes \underline{B}_1 \ \underline{A} \otimes \underline{B}_2 \ ... \underline{A} \otimes \underline{B}_n \right] \qquad (30)$$

Normally this matrix is not the same as the matrix $\underline{A} \otimes \underline{B}$. As easily seen they are, if $\underline{A}$ is a column $\underline{a}$:

$$\underline{a} * \underline{\underline{B}} = \underline{a} \otimes \underline{\underline{B}} \qquad (31)$$

PROPOSITION 3.4. *If in (30) the matrix $\underline{\underline{A}}$ is $m \times n$ and all the matrices*

$\underline{\underline{B}}_k, k = 1, 2, ..., n$ *and $\underline{\underline{C}}$ are $n \times n$, then for all $\underline{h} \in R^n$*

$$(\underline{\underline{A}} * \underline{\underline{C}})\underline{h}^2 = (\underline{\underline{A}}\underline{h}) \otimes (\underline{\underline{C}}\underline{h}) \quad and \quad (\underline{\underline{A}} * \underline{\underline{B}})\underline{h}^3 = (\underline{\underline{A}}\underline{h}) \otimes (\underline{\underline{B}}\underline{h}^2) \qquad (32)$$

PROPOSITION 3.5. *Let $\Omega \subseteq R^s$ be an open set and let $\underline{\underline{A}}(\underline{x}), \underline{\underline{B}}(\underline{x}), \underline{x} \in \Omega$ be $C^l$*

*matrix functions of form $m \times n$ and $p \times q$, respectively. Then*

$$\underline{D}(\underline{\underline{A}}\underline{\underline{B}}) = \underline{\underline{A}}(\underline{D}\underline{\underline{B}}) + \underline{D}\underline{\underline{A}}(\underline{\underline{E}}_s \otimes \underline{\underline{B}}) \quad if \quad n = p \qquad (33)$$

$$\underline{D}(\underline{\underline{A}} \otimes \underline{\underline{B}}) = \underline{\underline{A}} * \underline{D}\underline{\underline{B}} + (\underline{D}\underline{\underline{A}}) \otimes \underline{\underline{B}} \qquad (34)$$

*If in particular $\underline{\underline{A}}$ is an $m \times 1$ column $\underline{a}$ then (34) becomes*

$$\underline{D}(\underline{a} \otimes \underline{\underline{B}}) = \underline{a} \otimes \underline{D}\underline{\underline{B}} + (\underline{D}\underline{a}) \otimes \underline{\underline{B}} \qquad (35)$$

*and if m=1 then (35) becomes*

$$\underline{D}(a\underline{\underline{B}}) = a\underline{D}\underline{\underline{B}} + (\nabla a)^T \otimes \underline{\underline{B}} \qquad (36)$$

PROPOSITION 3.6. *Let $\Omega \subseteq R^s$ and $\Omega_1 \subseteq R^p$ be open sets and let, respectively,*

$\underline{\underline{A}}(\underline{x}), \underline{x} \in \Omega_1$ *and $\underline{\varphi} : \Omega \to \Omega_1$ be a $C^l$ $m \times n$ matrix function and a $C^l$ vector*

*function. Then*

$$\underline{D}(\underline{\underline{A}} \circ \underline{\varphi}) = [(\underline{D}\underline{\underline{A}}) \circ \underline{\varphi}][(\underline{D}\underline{\varphi}) \otimes \underline{\underline{E}}_n] \qquad (37)$$

COROLLARY 3.2. *If $\underline{f} : \Omega_1 \to R^n$ is a $C^q$ vector function then*

$$\underline{D}((\underline{D}^k \underline{f}) \circ \underline{\varphi}) = [(\underline{D}^{k+1} \underline{f}) \circ \underline{\varphi}][(\underline{D}\underline{\varphi}) \otimes \underline{\underline{E}}_{kp}] \quad , k = 1, 2, ..., q \qquad (38)$$

*Proof.* Put $\underline{\underline{A}} = \underline{D}^k \underline{f}$ in (37) and observe that n=kp and that $\underline{\underline{E}}_{kp} = \underline{\underline{E}}^{\otimes k}$.

PROPOSITION 3.7. *Let $\underline{\underline{C}} = [\underline{\underline{C}}_1 \underline{\underline{C}}_2 ... \underline{\underline{C}}_n]$ be an $n \times n^2$ matrix for which there*

*holds $\underline{\underline{C}}\underline{h}^2 = \alpha(\underline{h})\underline{h}$ for all $\underline{h} \in R^n$, where $\alpha(\underline{h}), \underline{h} \in R^n$ is some scalar*

*function. Then there exists an $n \times 1$ column $\underline{a}$ such that $\alpha(\underline{h}) = \underline{a}^T \underline{h}$ for all*

$\underline{h} \neq \underline{0}$. *If in particular $\underline{\underline{C}}$ is c-symmetric then $\underline{\underline{C}} = \frac{1}{2}(\underline{a}^T \otimes \underline{\underline{E}} + \underline{\underline{E}} \otimes \underline{a}^T)$.*

### 4. Proof of formula (9)

We consider the differential system

$$\frac{d\underline{x}}{dt} = \underline{f}(t, \underline{x}), (t, \underline{x}) \in I \times R^n \qquad (39)$$

where $I \subseteq R$ is an open interval and $\underline{f}(t,\underline{x}), (t,\underline{x}) \in I \times R^n$ is a continuous function, which is C³ in $\underline{x}$. Let $\underline{\varphi}(t,\tau,\underline{\xi}), t \in I(\tau,\underline{\xi})$ denote the maximal

solution of (42) through the point $(\tau,\underline{\xi}) \in I \times R^n$. It is well known that

$$\frac{\partial}{\partial t} \underline{D}\underline{\varphi}(t,\tau,\underline{\xi}) = \underline{D}\underline{f}(t,\underline{\varphi}(t,\tau,\underline{\xi}))\underline{D}\underline{\varphi}(t,\tau,\underline{\xi}), t \in I(\tau,\underline{\xi}) \tag{40}$$

see [5, p 40 ff]. Recall also that $\underline{d}\underline{\varphi}(t,\tau,\underline{\xi},\underline{h}) = \underline{D}\underline{\varphi}(t,\tau,\underline{\xi})\underline{h}$. Omitting all

arguments we have

$$\frac{\partial}{\partial t} \underline{D}\underline{\varphi} = \underline{D}\underline{f}\,\underline{D}\underline{\varphi} \tag{41}$$

$$\frac{\partial}{\partial t} \underline{d}\underline{\varphi} = \underline{D}\underline{f}\,\underline{d}\underline{\varphi} \tag{42}$$

In the following $\underline{D}_{\underline{\xi}}$ indicates differentiation with respect to $\underline{\xi}$, where there could be risk of misunderstanding (recall the $\underline{D}^k$ before $\underline{f}$ indicates differentiation with respect to $\underline{x}$ and before $\underline{\varphi}$ with respect to $\underline{\xi}$). We now derive the second order and third order analogues of (41) and (42).

Differentiating (41) with respect to $\underline{\xi}$ we get from (33) in Proposition 3.5

$$\frac{\partial}{\partial t} \underline{D}^2\underline{\varphi} = \underline{D}\underline{f}\,\underline{D}^2\underline{\varphi} + \underline{D}_{\underline{\xi}}(\underline{D}\underline{f})(\underline{E} \otimes \underline{D}\underline{\varphi})$$

Next we use Proposition 3.6 and the product rule obtaining

$$\frac{\partial}{\partial t} \underline{D}^2\underline{\varphi} = \underline{D}\underline{f}\,\underline{D}^2\underline{\varphi} + \underline{D}^2\underline{f}(\underline{D}\underline{\varphi} \otimes \underline{E})(\underline{E} \otimes \underline{D}\underline{\varphi}) = \underline{D}\underline{f}\,\underline{D}^2\underline{\varphi} + \underline{D}^2\underline{f}(\underline{D}\underline{\varphi} \otimes \underline{D}\underline{\varphi}) \tag{43}$$

where $\underline{D}^2\underline{f}$ is taken at $(t,\underline{\varphi}(t,\tau,\underline{\xi})$. From the product rule it follows that

$(\underline{D}\underline{\varphi} \otimes \underline{D}\underline{\varphi})\underline{h}^2 = (\underline{D}\underline{\varphi}\underline{h}) \otimes (\underline{D}\underline{\varphi}\underline{h}) = (\underline{d}\underline{\varphi})^2$. Hence (43) implies that

$$\frac{\partial}{\partial t} \underline{d}^2\underline{\varphi} = \underline{D}\underline{f}\,\underline{d}^2\underline{\varphi} + \underline{D}^2\underline{f}(\underline{d}\underline{\varphi})^2 \tag{44}$$

Next the same procedure is repeated with (43) instead of (41). Also using

Corollary 3.2 we get

$$\frac{\partial}{\partial t} \underline{D}^3\underline{\varphi} = \underline{D}_{\underline{\xi}}(\underline{D}\underline{f}\,\underline{D}^2\underline{\varphi}) + \underline{D}_{\underline{\xi}}\left((\underline{D}^2\underline{f})(\underline{D}\underline{\varphi} \otimes \underline{D}\underline{\varphi})\right) = \underline{D}\underline{f}\,\underline{D}^3\underline{\varphi} + (\underline{D}^2\underline{f})(\underline{D}\underline{\varphi} \otimes \underline{E})(\underline{E} \otimes \underline{D}^2\underline{\varphi})$$

$$+ \underline{D}^2\underline{f}\,\underline{D}_{\underline{\xi}}(\underline{D}\underline{\varphi} \otimes \underline{D}\underline{\varphi}) + \underline{D}^3\underline{f}(\underline{D}\underline{\varphi} \otimes \underline{E}_{2n})(\underline{E} \otimes \underline{D}\underline{\varphi} \otimes \underline{D}\underline{\varphi})$$

By the product rule we have $(\underline{D}\underline{\varphi} \otimes \underline{E})(\underline{E} \otimes \underline{D}^2\underline{\varphi}) = \underline{D}\underline{\varphi} \otimes \underline{D}^2\underline{\varphi}$ and in the

same way

$(\underline{D}\underline{\varphi} \otimes \underline{\underline{E}}_{2n})(\underline{\underline{E}} \otimes \underline{D}\underline{\varphi} \otimes \underline{D}\underline{\varphi}) = (\underline{D}\underline{\varphi} \otimes \underline{\underline{E}} \otimes \underline{\underline{E}})(\underline{\underline{E}} \otimes \underline{D}\underline{\varphi} \otimes \underline{D}\underline{\varphi}) = \underline{D}\underline{\varphi} \otimes \underline{D}\underline{\varphi} \otimes \underline{D}\underline{\varphi}$

Inserting in the above expression we get

$$\frac{\partial}{\partial t}\underline{D}^3\underline{\varphi} = \underline{D}\underline{f}\,\underline{D}^3\underline{\varphi} + \underline{D}^2\underline{f}(\underline{D}\underline{\varphi} \otimes \underline{D}^2\underline{\varphi}) + \underline{D}^2\underline{f}(\underline{D}\underline{\varphi} * \underline{D}^2\underline{\varphi} + \underline{D}^2\underline{\varphi} \otimes \underline{D}\underline{\varphi}) \qquad (45)$$
$$+ \underline{D}^3\underline{f}(\underline{D}\underline{\varphi} \otimes \underline{D}\underline{\varphi} \otimes \underline{D}\underline{\varphi})$$

Using (32) and next (24) we finally get

$$\frac{\partial}{\partial t}\underline{d}^3\underline{\varphi} = \underline{D}\underline{f}\,\underline{d}^3\underline{\varphi} + 3\underline{D}^2\underline{f}(\underline{d}^2\underline{\varphi} \otimes \underline{d}\underline{\varphi}) + \underline{D}^3\underline{f}(\underline{d}\underline{\varphi})^3 \qquad (46)$$

We are now in a position to prove formula (9). The proof will be easier if the factor $\underline{D}\underline{\varphi} \otimes \underline{D}\underline{\varphi}$ at the right hand side of (9) is moved to the left hand side.

To this end we introduce the following abbreviations $\underline{\underline{\Psi}} = \left(\underline{D}\underline{\varphi}\right)^{-1}$,

$\underline{P} = \underline{\underline{\Psi}}\underline{D}^2\underline{\varphi}$ and $\underline{Q} = \underline{\underline{\Psi}}\underline{D}^3\underline{\varphi}$, whence

$$\underline{P}\underline{h}^2 = \underline{\underline{\Psi}}\underline{d}^2\underline{\varphi} \quad \text{and} \quad \underline{Q}\underline{h}^3 = \underline{\underline{\Psi}}\underline{d}^3\underline{\varphi} \qquad (47)$$

Note that $\underline{d}\underline{\varphi} = \underline{D}\underline{\varphi}\underline{h}$ implies that $\underline{h} = \underline{\underline{\Psi}}\underline{d}\underline{\varphi}$. Moreover, (41) implies that

$$\frac{\partial}{\partial t}\underline{\underline{\Psi}} = -\left(\underline{D}\underline{\varphi}\right)^{-1}\frac{\partial}{\partial t}\underline{D}\underline{\varphi}\left(\underline{D}\underline{\varphi}\right)^{-1} = -\underline{\underline{\Psi}}\underline{D}\underline{f} \qquad (48)$$

If we multiply with $(\underline{D}\underline{\varphi} \otimes \underline{D}\underline{\varphi})^{-1} = (\underline{D}\underline{\varphi})^{-1} \otimes (\underline{D}\underline{\varphi})^{-1} = \underline{\underline{\Psi}} \otimes \underline{\underline{\Psi}}$ on both sides of (9) the right hand side becomes a sum of integrals. By the product rule the left hand side becomes

$$(\underline{\underline{\Psi}} \otimes \underline{\underline{\Psi}})\left[(\underline{D}^3\underline{\varphi}\underline{h}^3) \otimes \underline{d}\underline{\varphi} + \underline{d}\underline{\varphi} \otimes (\underline{D}^3\underline{\varphi}\underline{h}^3) - 3(\underline{d}^2\underline{\varphi}) \otimes (\underline{d}^2\underline{\varphi})\right]$$
$$= \underline{Q}\underline{h}^3 \otimes \underline{h} + \underline{h} \otimes \underline{Q}\underline{h}^3 - 3\underline{P}\underline{h}^2 \otimes \underline{P}\underline{h}^2$$

We show that the time derivative of this expression equals the sum of the integrands in (9), and then formula (9) will be proved. First we get from (44)

$$\frac{\partial}{\partial t}\left(\underline{P}\underline{h}^2\right) = \frac{\partial}{\partial t}\left(\underline{\underline{\Psi}}\underline{d}^2\underline{\varphi}\right) = -\underline{\underline{\Psi}}\underline{D}\underline{f}\,\underline{d}^2\underline{\varphi} + \underline{\underline{\Psi}}\underline{D}\underline{f}\,\underline{d}^2\underline{\varphi} + \underline{\underline{\Psi}}\underline{D}^2\underline{f}(\underline{d}\underline{\varphi})^2 = \underline{\underline{\Psi}}\underline{D}^2\underline{f}(\underline{d}\underline{\varphi})^2 \quad (49)$$

and next from (46)

$$\frac{\partial}{\partial t}\left(\underline{Q}\underline{h}^3\right) = \frac{\partial}{\partial t}\left(\underline{\underline{\Psi}}\underline{d}^3\underline{\varphi}\right) = -\underline{\underline{\Psi}}\underline{D}\underline{f}\,\underline{d}^3\underline{\varphi} + \underline{\underline{\Psi}}(\underline{D}\underline{f}\,\underline{d}^3\underline{\varphi} + 3\underline{D}^2\underline{f}(\underline{d}^2\underline{\varphi} \otimes \underline{d}\underline{\varphi}) + \underline{D}^3\underline{f}(\underline{d}\underline{\varphi})^3)$$
$$= \underline{\underline{\Psi}}(3\underline{D}^2\underline{f}(\underline{d}^2\underline{\varphi} \otimes \underline{d}\underline{\varphi}) + \underline{D}^3\underline{f}(\underline{d}\underline{\varphi})^3)$$

Using the product rule twice we get

$$\frac{\partial}{\partial t}\left(\underline{\underline{Q}}h^3 \otimes \underline{h}\right) = \frac{\partial}{\partial t}\left(\underline{\underline{Q}}h^3\right) \otimes \underline{h} = \left(\underline{\Psi}(3\underline{\underline{D}}^2\underline{f}(\underline{d}^2\varphi \otimes \underline{d\varphi}) + \underline{D}^3\underline{f}(\underline{d\varphi})^3)\right) \otimes \underline{\Psi}\underline{d\varphi}$$

$$= \left(\underline{\Psi} \otimes \underline{\underline{\Psi}}\right)\left((3\underline{\underline{D}}^2\underline{f}(\underline{d}^2\varphi \otimes \underline{d\varphi}) + \underline{D}^3\underline{f}(\underline{d\varphi})^3)) \otimes \underline{E}\underline{d\varphi}\right)$$

$$= \left(\underline{\Psi} \otimes \underline{\underline{\Psi}}\right)\left((3(\underline{\underline{D}}^2\underline{f} \otimes \underline{E})(\underline{d}^2\varphi \otimes (\underline{d\varphi})^2) + (\underline{D}^3\underline{f} \otimes \underline{E})(\underline{d\varphi})^4\right) \quad (50)$$

and analogously

$$\frac{\partial}{\partial t}\left(\underline{h} \otimes \underline{\underline{Q}}h^3\right) = \left(\underline{\Psi} \otimes \underline{\underline{\Psi}}\right)\left((3(\underline{E} \otimes \underline{\underline{D}}^2\underline{f})(\underline{d\varphi} \otimes \underline{d}^2\varphi \otimes \underline{d\varphi}) + (\underline{E} \otimes \underline{D}^3\underline{f})(\underline{d\varphi})^4\right)$$

By (26) we have $(\underline{E} \otimes \underline{\underline{D}}^2\underline{f})(\underline{d\varphi} \otimes \underline{d}^2\varphi \otimes \underline{d\varphi}) = (\underline{E} \otimes \underline{\underline{D}}^2\underline{f})((\underline{d\varphi})^2 \otimes \underline{d}^2\varphi)$

since $\underline{\underline{D}}^2\underline{f}$ is c-symmetric. This rearrangement is crucial . We get

$$\frac{\partial}{\partial t}\left(\underline{h} \otimes \underline{\underline{Q}}h^3\right) = \left(\underline{\Psi} \otimes \underline{\underline{\Psi}}\right)\left((3(\underline{E} \otimes \underline{\underline{D}}^2\underline{f})((\underline{d\varphi})^2 \otimes \underline{d}^2\varphi) + (\underline{E} \otimes \underline{D}^3\underline{f})(\underline{d\varphi})^4\right) \quad (51)$$

We further get

$$\frac{\partial}{\partial t}\left((\underline{P}h^2) \otimes (\underline{P}h^2)\right) = (\underline{\Psi}\underline{d}^2\varphi) \otimes (\underline{\Psi}\underline{D}^2\underline{f}(\underline{d\varphi})^2) + (\underline{\Psi}\underline{D}^2\underline{f}(\underline{d\varphi})^2) \otimes (\underline{\Psi}\underline{d}^2\varphi)$$

$$= \left(\underline{\Psi} \otimes \underline{\underline{\Psi}}\right)\left((\underline{E}\underline{d}^2\varphi) \otimes \underline{\underline{D}}^2\underline{f}(\underline{d\varphi})^2 + (\underline{\underline{D}}^2\underline{f}(\underline{d\varphi})^2) \otimes (\underline{E}\underline{d}^2\varphi)\right)$$

$$= \left(\underline{\Psi} \otimes \underline{\underline{\Psi}}\right)\left((\underline{E} \otimes \underline{\underline{D}}^2\underline{f})(\underline{d}^2\varphi \otimes (\underline{d\varphi})^2) + (\underline{\underline{D}}^2\underline{f} \otimes \underline{E})((\underline{d\varphi})^2 \otimes \underline{d}^2\varphi)\right)$$

Combining this and (50) and (51) we get

$$\frac{\partial}{\partial t}\left((\underline{\underline{Q}}h^3) \otimes \underline{h} + \underline{h} \otimes (\underline{\underline{Q}}h^3) - 3(\underline{P}h^2) \otimes (\underline{P}h^2)\right)$$

$$= \left(\underline{\Psi} \otimes \underline{\underline{\Psi}}\right)[3(\underline{\underline{D}}^2\underline{f} \otimes \underline{E})(\underline{d}^2\varphi \otimes (\underline{d\varphi})^2) + (\underline{D}^3\underline{f} \otimes \underline{E})(\underline{d\varphi})^4 + 3(\underline{E} \otimes \underline{\underline{D}}^2\underline{f})((\underline{d\varphi})^2 \otimes \underline{d}^2\varphi)$$

$$+ (\underline{E} \otimes \underline{D}^3\underline{f})(\underline{d\varphi})^4 - 3(\underline{E} \otimes \underline{\underline{D}}^2\underline{f})(\underline{d}^2\varphi \otimes (\underline{d\varphi})^2) - 3(\underline{\underline{D}}^2\underline{f} \otimes \underline{E})((\underline{d\varphi})^2 \otimes \underline{d}^2\varphi)]$$

$$= \left(\underline{\Psi} \otimes \underline{\underline{\Psi}}\right)[(\underline{D}^3\underline{f} \otimes \underline{E} + \underline{E} \otimes \underline{D}^3\underline{f})(\underline{d\varphi})^4$$

$$+ 3(\underline{\underline{D}}^2\underline{f} \otimes \underline{E} - \underline{E} \otimes \underline{\underline{D}}^2\underline{f})(\underline{d}^2\varphi \otimes (\underline{d\varphi})^2 - (\underline{d\varphi})^2 \otimes \underline{d}^2\varphi)] \quad (52)$$

Note that $t = \tau$ gives $\underline{\underline{D}}\varphi(\tau, \tau, \underline{\xi}) = \underline{\Phi}(\tau, \tau, \underline{\xi}) = \underline{E}$, whence $\underline{\underline{D}}^2\varphi(\tau, \tau, \underline{\xi})$ and

$\underline{\underline{D}}^2\varphi(\tau, \tau, \underline{\xi})$ are both zero matrices. Hence also $\underline{P}$ and $\underline{\underline{Q}}$ are zero matrices for

$t = \tau$. Integrating the left hand side of the first equality sign in (52) we get

$$\int_\tau^t \frac{\partial}{\partial t}\left((\underline{\underline{Q}}h^3) \otimes \underline{h} + \underline{h} \otimes (\underline{\underline{Q}}h^3) - 3(\underline{P}h^2) \otimes (\underline{P}h^2)\right) ds$$

$$= (\underline{\underline{Q}}h^3) \otimes \underline{h} + \underline{h} \otimes (\underline{\underline{Q}}h^3) - 3(\underline{P}h^2) \otimes (\underline{P}h^2)$$

$$= \left(\underline{\Psi} \otimes \underline{\underline{\Psi}}\right)\left(\underline{d}^3\varphi \otimes \underline{d\varphi} + \underline{d\varphi} \otimes \underline{d}^3\varphi - 3(\underline{d}^2\varphi)^2\right) \quad (53)$$

Integrating the right hand side of the last equality sign in (52) and identifying this with (53), we get an expression which after insertion of $\underline{\Psi} = \left(\underline{\underline{D}}\varphi\right)^{-1}$ and

multiplication with $\underline{\underline{D}}\varphi \otimes \underline{\underline{D}}\varphi$ becomes

$$\underline{d}^3\underline{\varphi}\otimes\underline{d\varphi}+\underline{d\varphi}\otimes\underline{d}^3\underline{\varphi}-3\left[\underline{d}^2\underline{\varphi}\right]^2=\underline{\underline{D}\varphi}\otimes\underline{\underline{D}\varphi}\left\{\int_\tau^t\left[\underline{\underline{D}\varphi}\times\underline{\underline{D}\varphi}\right]^{-1}\left(\underline{\underline{D}}^3\underline{f}\otimes\underline{\underline{E}}+\underline{\underline{E}}\otimes\underline{\underline{D}}^3\underline{f}\right)\!\left(\underline{d\varphi}\right)^4ds\right\}$$

$$+3\left\{\int_\tau^t\left[\underline{\underline{D}\varphi}\otimes\underline{\underline{D}\varphi}\right]^{-1}\left(\underline{\underline{D}}^2\underline{f}\otimes\underline{\underline{E}}-\underline{\underline{E}}\otimes\underline{\underline{D}}^2\underline{f}\right)\!\left(\underline{d}^2\underline{\varphi}\otimes\left(\underline{d\varphi}\right)^2-\left(\underline{d\varphi}\right)^2\otimes\underline{d}^2\underline{\varphi}\right)ds\right\}$$

which is formula (9). In the integrands $\underline{\underline{D}}^k\underline{f}$ shall be taken in the point

$(s,\underline{\varphi}(s,\tau,\underline{\xi}))$, and $\underline{\underline{D}}^k\underline{\varphi}$ shall be taken in the point $(s,\tau,\underline{\xi})$ in the integrands and in the point $(t,\tau,\underline{\xi})$ outside the integrands. Recall also that $\underline{d}^k\underline{\varphi}=\underline{\underline{D}}^k\underline{\varphi}\underline{h}$. This completes the proof.

It should be noted that (49) in the same way gives

$$\underline{\underline{P}}\underline{h}^2=\underline{\underline{\Psi}}\underline{\underline{D}}^2\underline{\varphi}\underline{h}^2=\int_\tau^t\underline{\underline{\Psi}}\underline{\underline{D}}^2\underline{f}\left(\underline{d\varphi}\right)^2ds=\int_\tau^t\underline{\underline{\Psi}}\underline{\underline{D}}^2\underline{f}\left(\underline{\underline{D}\varphi}\otimes\underline{\underline{D}\varphi}\right)ds\,\underline{h}^2\qquad(54)$$

It follows from Corollary 2.2 and the remarks thereafter that $\underline{h}^2$ may be omitted in (54). Inserting $\underline{\underline{\Psi}}=\left(\underline{\underline{D}\varphi}\right)^{-1}$ we the obtain

$$\left(\underline{\underline{D}\varphi}\right)^{-1}\underline{\underline{D}}^2\underline{\varphi}=\int_\tau^t\left(\underline{\underline{D}\varphi}\right)^{-1}\underline{\underline{D}}^2\underline{f}\left(\underline{\underline{D}\varphi}\otimes\underline{\underline{D}\varphi}\right)ds$$

which after multiplication with $\underline{\underline{D}\varphi}$ on both sides becomes formula (8)

(recall that $\underline{\underline{D}\varphi}=\underline{\underline{\Phi}}$). By the product rule formula (9) can be written

$$\left(\underline{\underline{D}}^3\underline{\varphi}\otimes\underline{\underline{D}\varphi}+\underline{\underline{D}\varphi}\otimes\underline{\underline{D}}^3\underline{\varphi}-3\underline{\underline{D}}^2\underline{\varphi}\otimes\underline{\underline{D}}^2\underline{\varphi}\right)\underline{h}^4=$$

$$\underline{\underline{D}\varphi}\otimes\underline{\underline{D}\varphi}\left\{\int_\tau^t\left[\underline{\underline{D}\varphi}\otimes\underline{\underline{D}\varphi}\right]^{-1}\left(\underline{\underline{D}}^3\underline{f}\otimes\underline{\underline{E}}+\underline{\underline{E}}\otimes\underline{\underline{D}}^3\underline{f}\right)\!\left(\underline{\underline{D}\varphi}\right)^{\otimes4}ds\right\}$$

$$+3\left\{\int_\tau^t\left[\underline{\underline{D}\varphi}\otimes\underline{\underline{D}\varphi}\right]^{-1}\left(\underline{\underline{D}}^2\underline{f}\otimes\underline{\underline{E}}-\underline{\underline{E}}\otimes\underline{\underline{D}}^2\underline{f}\right)\!\left(\underline{\underline{D}}^2\underline{\varphi}\otimes\underline{\underline{D}\varphi}\otimes\underline{\underline{D}\varphi}-\underline{\underline{D}\varphi}\otimes\underline{\underline{D}\varphi}\otimes\underline{d}^2\underline{\varphi}\right)ds\right\}\underline{h}^4$$

Contrary to (54) the factor $\underline{h}^4$ cannot be omitted in this expression. This does

not mean that one cannot derive a formula like (8) without the factor $\underline{h}^4$. In

fact in the same way as above we could calculate an expression for the time

derivative of

$$\left(\underline{\underline{D}\varphi}\otimes\underline{\underline{D}\varphi}\right)^{-1}\left(\underline{\underline{D}}^3\underline{\varphi}\otimes\underline{\underline{D}\varphi}+\underline{\underline{D}\varphi}\otimes\underline{\underline{D}}^3\underline{\varphi}-3\underline{\underline{D}}^2\underline{\varphi}\otimes\underline{\underline{D}}^2\underline{\varphi}\right)$$

$$=\left[\left(\underline{\underline{D}\varphi}\right)^{-1}\underline{\underline{D}}^3\underline{\varphi}\right]\otimes\underline{\underline{E}}+\underline{\underline{E}}\otimes\left[\left(\underline{\underline{D}\varphi}\right)^{-1}\underline{\underline{D}}^3\underline{\varphi}\right]-3\left[\left(\underline{\underline{D}\varphi}\right)^{-1}\underline{\underline{D}}^2\underline{\varphi}\right]\otimes\left[\left(\underline{\underline{D}\varphi}\right)^{-1}\underline{\underline{D}}^2\underline{\varphi}\right]$$

However, this expression will be extensive since the time derivative of e.g.

$\left[\left(\underline{\underline{D}\varphi}\right)^{-1}\underline{\underline{D}}^3\underline{\varphi}\right]\otimes\underline{\underline{E}}$ will contain more terms than the time derivative of the corresponding term $\left(\underline{\underline{Q}}\underline{h}^3\right)\otimes\underline{h}$ in the calculations leading to (52).

For later use we note that (54) implies that

$$\underline{d}^2\underline{\varphi} = \underline{\underline{D\varphi}} \int_\tau^t \left(\underline{\underline{D\varphi}}\right)^{-1} \underline{\underline{D}}^2 \underline{f}(\underline{d\varphi})^2 ds \qquad (55)$$

### 5. The vector Riccati equation

As mentioned in the introduction a vector Riccati equation has the form

$$\frac{d\underline{x}}{dt} = \underline{a}(t) + \underline{\underline{B}}(t)\underline{x} + \left(\underline{\underline{c}}^T(t)\underline{x}\right)\underline{x}, \; t \in I \qquad (13)$$

where $\underline{a}(t), t \in I$, $\underline{c}(t), t \in I$ are continuous vector functions and $\underline{\underline{B}}(t), t \in I$ is a continuous matrix function.

THEOREM 5.1. *System (39) is a vector Riccati equation if and only if*

$\underline{d}^3\underline{\varphi} \otimes \underline{d\varphi} + \underline{d\varphi} \otimes \underline{d}^3\underline{\varphi} - 3(\underline{d}^2\underline{\varphi})^2 = \underline{0}$ *for all* $t, \tau \in I$, $\underline{\xi} \in R^n$ *and all* $\underline{h} \in R^n$.

*Proof.* Suppose first that $\underline{d}^3\underline{\varphi} \otimes \underline{d\varphi} + \underline{d\varphi} \otimes \underline{d}^3\underline{\varphi} = 3(\underline{d}^2\underline{\varphi})^2$ holds identically. From (28) in Proposition 3.3 follows that

$$\underline{d}^2\underline{\varphi} = \alpha(t, \tau, \underline{\xi}, \underline{h})\underline{d\varphi} \qquad (56)$$

where $\alpha(t, \tau, \underline{\xi}, \underline{h})$ is a scalar function. Consequently, $\underline{d\varphi}$ and $\underline{d}^2\underline{\varphi}$ Kronecker commute. As $\underline{\underline{D\varphi}} \otimes \underline{\underline{D\varphi}}$ is nonsingular we derive from formula (9) that

$$\int_\tau^t \left[\underline{\underline{D\varphi}} \otimes \underline{\underline{D\varphi}}\right]^{-1} \left(\underline{\underline{D}}^3 \underline{f} \otimes \underline{\underline{E}} + \underline{\underline{E}} \otimes \underline{\underline{D}}^3 \underline{f}\right)(\underline{d\varphi})^4 ds = \underline{0} \qquad (57)$$

identically. We infer that the integrand in (57) is identically $\underline{0}$, whence

$$\left(\underline{\underline{D}}^3 \underline{f} \otimes \underline{\underline{E}} + \underline{\underline{E}} \otimes \underline{\underline{D}}^3 \underline{f}\right)\!\left(\underline{d\varphi}\right)^4 = \left(\underline{\underline{D}}^3 \underline{f}(\underline{d\varphi})^3\right) \otimes \underline{d\varphi} + \underline{d\varphi} \otimes \left(\underline{\underline{D}}^3 \underline{f}(\underline{d\varphi})^3\right) = \underline{0}$$

identically. Since $\underline{d\varphi} = \underline{\underline{D\varphi}}\underline{h} \neq \underline{0}$ if $\underline{h} \neq \underline{0}$ we get from (27) in Proposition 3.3 that $\underline{\underline{D}}^3 \underline{f}(\underline{d\varphi})^3 = \underline{\underline{D}}^3 \underline{f}(\underline{\underline{D\varphi}}\underline{h})^3 = \underline{0}$ identically. Trivially, this is also true if $\underline{h} = \underline{0}$. By Corollary 2.1 we infer that $\underline{\underline{D}}^3 \underline{f}(t, \underline{\varphi}(t, \tau, \underline{\xi})) = \underline{0}$ identically. Hence $\underline{\underline{D}}^3 \underline{f}(t, \underline{x}) = \underline{0}$ for all $(t, \underline{x}) \in I \times R^n$. Using Taylors formula (18) we may write

$$\underline{f}(t, \underline{x}) = \underline{a}(t) + \underline{\underline{B}}(t)\underline{x} + \frac{1}{2}\underline{\underline{C}}(t)\underline{x}^2 , (t, \underline{x}) \in I \times R^n \qquad (58)$$

where $\underline{a}(t) = \underline{f}(t, \underline{0})$, $\underline{\underline{B}}(t) = \underline{\underline{D}}\underline{f}(t, \underline{0})$ and $\underline{\underline{C}}(t) = \underline{\underline{D}}^2 \underline{f}(t, \underline{0})$. We may write (56) as

$$\left(\underline{\underline{D\varphi}}\right)^{-1} \underline{\underline{D}}^2 \underline{\varphi}\underline{h}^2 = \alpha(t, \tau, \underline{\xi}, \underline{h})\underline{h} \qquad (59)$$

From Proposition 3.7 follows that $\alpha(t,\tau,\underline{\xi},\underline{h}) = \underline{a}^T(t,\tau,\underline{\xi})\underline{h}$ for some $n \times 1$ column $\underline{a}(t,\tau,\underline{\xi})$. Combining (55) , (56) and (59) we get

$$\alpha(t,\tau,\underline{\xi},\underline{h})d\underline{\varphi} = \underline{\underline{D}}\underline{\varphi}\big((\underline{a}^T(t,\tau,\underline{\xi})\underline{h})\underline{h}\big) = \underline{\underline{D}}\underline{\varphi}\int_\tau^t \big(\underline{\underline{D}}\underline{\varphi}\big)^{-1}\underline{\underline{C}}(s)\big(\underline{\underline{D}}\underline{\varphi} \otimes \underline{\underline{D}}\underline{\varphi}\big)\,ds\,\underline{h}^2$$

since $\underline{\underline{D}}^2\underline{f}(t,\underline{x}) = \underline{\underline{C}}(t)$ by (58). We infer that

$$(\underline{a}^T(t,\tau,\underline{\xi})\underline{h})\underline{h} = \int_\tau^t \big(\underline{\underline{D}}\underline{\varphi}\big)^{-1}\underline{\underline{C}}(s)\big(\underline{\underline{D}}\underline{\varphi} \otimes \underline{\underline{D}}\underline{\varphi}\big)\,ds\,\underline{h}^2 \qquad (60)$$

It is obvious that (60) implies that $\underline{a}(t,\tau,\underline{\xi})$ is differentiable with respect to $t$, and we get

$$(\frac{\partial}{\partial t}\underline{a}^T(t,\tau,\underline{\xi})\underline{h})\underline{h} = \big(\underline{\underline{D}}\underline{\varphi}(t,\tau,\underline{\xi})\big)^{-1}\underline{\underline{C}}(t)\big(\underline{\underline{D}}\underline{\varphi}(t,\tau,\underline{\xi}) \otimes \underline{\underline{D}}\underline{\varphi}(t,\tau,\underline{\xi})\big)\underline{h}^2 \qquad (61)$$

There holds $\underline{\underline{D}}\underline{\varphi}(t,t,\underline{\xi}) = \underline{\underline{E}}$ since $\underline{\varphi}(t,t,\underline{\xi}) = \underline{\xi}$ for all $(t,\underline{\xi}) \in I \times R^n$. Putting $t = \tau$ in (61) we therefore get

$$\left(\left(\frac{\partial}{\partial t}\underline{a}^T(t,\tau,\underline{\xi})\right)_{t=\tau}\underline{h}\right)\underline{h} = \underline{\underline{C}}(\tau)\underline{h}^2 \qquad (62)$$

Choosing $\underline{h} = \underline{e}_i$ in (62) and multiplying both sides from the left with $\underline{e}_i^T$ we get

$$\left(\frac{\partial}{\partial t}a_i(t,\tau,\underline{\xi})\right)_{t=\tau} = \underline{e}_i^T\left(\frac{\partial}{\partial t}a_i(t,\tau,\underline{\xi})\right)_{t=\tau}\underline{e}_i = \underline{e}_i^T\underline{\underline{C}}(\tau)\underline{e}_i^2 = 2c_i(\tau)$$

where the last equality sign indicates the definition of $c_i(\tau)$. Note that this is valid independently of $\underline{\xi}$. Hence (62) may be written

$$\big(2\underline{c}^T(\tau)\underline{h}\big)\underline{h} = \underline{\underline{C}}(\tau)\underline{h}^2 \qquad (63)$$

for all $\tau \in I$, $\underline{h} \in R^n$. Substituting $\tau = t$, $\underline{h} = \underline{x}$ in (63) and next inserting in (58) we obtain the desired expression for $\underline{f}(t,\underline{x})$.

Suppose next that $\underline{f}(t,\underline{x}) = \underline{a}(t) + \underline{\underline{B}}(t)\underline{x} + \big(\underline{c}^T(t)\underline{x}\big)\underline{x}$, $(t,\underline{x}) \in I \times R^n$. Then $\underline{\underline{D}}^2\underline{f}(t,\underline{x}) = \underline{c}^T(t) \otimes \underline{\underline{E}} + \underline{\underline{E}} \otimes \underline{c}^T(t)$ and $\underline{\underline{D}}^3\underline{f}(t,\underline{x}) = \underline{\underline{0}}$. Hence, if $\underline{y}$ is a column then by the product rule

$$\underline{\underline{D}}^2\underline{f}(t,\underline{x})\underline{y}^2 = \big(\underline{c}^T(t) \otimes \underline{\underline{E}} + \underline{\underline{E}} \otimes \underline{c}^T(t)\big)\big(\underline{y} \otimes \underline{y}\big)$$
$$= \big(\underline{c}^T(t)\underline{y}\big) \otimes \underline{y} + \underline{y} \otimes \big(\underline{c}^T(t)\underline{b}\big) = 2\big(\underline{c}^T(t)\underline{y}\big)\underline{y}$$

Inserting the expression for $\underline{f}(t,\underline{x})$ in (55) and using this we get

$$\underline{d}^2\underline{\varphi} = \underline{\underline{D}}\underline{\varphi}\int_\tau^t\left(\underline{\underline{D}}\underline{\varphi}\right)^{-1}\left(2\underline{c}^T(s)\underline{d}\underline{\varphi}\right)\underline{\underline{D}}\underline{\varphi}\underline{h}ds = \underline{\underline{D}}\underline{\varphi}\int_\tau^t\left(2\underline{c}^T(s)\underline{d}\underline{\varphi}\right)\underline{h}\,ds = \left(\int_\tau^t 2\underline{c}^T(s)\underline{d}\underline{\varphi}\,ds\right)\underline{d}\underline{\varphi}$$

This implies that $\underline{d}^2\underline{\varphi}$ and $\underline{d}\underline{\varphi}$ Kronecker commute. Consequently, the right hand side of formula (9) becomes $\underline{0}$ identically. This proves the theorem.

### 6. Another characterization of the vector Riccati equation

By a fractional linear vector function of the variable $\underline{x} \in R^n$ we mean a vector function of the form

$$\underline{g}(\underline{x}) = \frac{\underline{\underline{A}}\underline{x} + \underline{\beta}}{\underline{\gamma}^T\underline{x} + \delta} \ , \ \underline{\gamma}^T\underline{x} + \delta \neq \underline{0} \tag{64}$$

where $\underline{\underline{A}}$ is an $n \times n$ matrix, $\underline{\beta}$ and $\underline{\gamma}$ are columns and $\delta$ a scalar. We assume that $\left[\underline{\gamma}^T\ \delta\right] \neq \underline{0}^T$.

PROPOSITION 6.1. *The general solution $\underline{\varphi}(t,\tau,\underline{\xi})$ of the vector Riccati equation (13) is a fractional linear function of $\underline{\xi}$, i.e., it has the form*

$$\underline{\varphi}(t,\tau,\underline{\xi}) = \frac{\underline{\underline{A}}(t,\tau)\underline{\xi} + \underline{\beta}(t,\tau)}{\underline{\gamma}^T(t,\tau)\underline{\xi} + \delta(t,\tau)}, t \in I(\tau,\underline{\xi}) \tag{65}$$

*Proof.* We introduce the linear differential system

$$\left|\begin{array}{c}\frac{d}{dt}\underline{y} \\ z\end{array}\right| = \left|\begin{array}{cc}\underline{\underline{B}}(t) & \underline{a}(t) \\ -\underline{c}^T(t) & 0\end{array}\right|\left|\begin{array}{c}\underline{y} \\ z\end{array}\right|, t \in I \tag{66}$$

Let

$$\underline{\underline{\Phi}}(t,\tau) = \left|\begin{array}{cc}\underline{\underline{\Psi}}(t,\tau) & \underline{\beta}(t,\tau) \\ \underline{\gamma}^T(t,\tau) & \delta(t,\tau)\end{array}\right|, t,\tau \in I \tag{67}$$

denote the fundamental matrix of system (66) for which $\underline{\underline{\Phi}}(t,\tau) = \underline{\underline{E}}_{n+1}$. Then the solution of system (66) through the point $(t,\underline{y},z) = (\tau,\underline{\xi},1)$ can be written

$$\left|\begin{array}{c}\underline{\psi}(t,\tau,\underline{\xi}) \\ \rho(t,\tau,\underline{\xi})\end{array}\right| = \underline{\underline{\Phi}}(t,\tau)\left[\begin{array}{c}\underline{\xi} \\ 1\end{array}\right] = \left|\begin{array}{c}\underline{\underline{\Psi}}(t,\tau)\underline{\xi} + \underline{\beta}(t,\tau) \\ \underline{\gamma}^T(t,\tau)\underline{\xi} + \delta(t,\tau)\end{array}\right|, t \in I$$

Since $\rho(\tau,\tau,\underline{\xi}) = 1$ there holds $\rho(t,\tau,\underline{\xi}) \neq 0$ for $t$ in some open interval around $t = \tau$. Let $J(\tau,\underline{\xi})$ denote the maximal subinterval of I for which $\rho(t,\tau,\underline{\xi}) \neq 0$. We put

$$\underline{x}(t,\tau,\underline{\xi}) = \frac{\underline{\psi}(t,\tau,\underline{\xi})}{\rho(t,\tau,\underline{\xi})}, \, t \in J(\tau,\underline{\xi})$$

Omitting all arguments we get

$$\frac{d\underline{x}}{dt} = -\frac{1}{\rho^2}\frac{d\rho}{dt}\underline{\psi} + \frac{1}{\rho}\frac{d\underline{\psi}}{dt} = -\frac{1}{\rho}(-\underline{c}^T\underline{\psi})\underline{x} + \frac{1}{\rho}\left(\underline{\underline{B}}\underline{\psi} + \rho\underline{a}\right) = (\underline{c}^T\underline{x})\underline{x} + \underline{\underline{B}}\underline{x} + \underline{a}, \, t \in J(\tau,\underline{\xi})$$

In other words, $\underline{x}(t,\tau,\underline{\xi})$ is a solution of the vector Riccati equation (13). Moreover,

$$\underline{x}(\tau,\tau,\underline{\xi}) = \frac{\underline{\psi}(\tau,\tau,\underline{\xi})}{\rho(\tau,\tau,\underline{\xi})} = \underline{\xi}$$

since $\underline{\underline{\Phi}}(t,\tau) = \underline{\underline{E}}_{n+1}$. Consequently, $\underline{x}(t,\tau,\underline{\xi}) = \underline{\varphi}(t,\tau,\underline{\xi})$ for all $t \in J(\tau,\underline{\xi})$. By a standard argument involving the Existence and Uniqueness Theorem we infer that $J(\tau,\underline{\xi}) = I(\tau,\underline{\xi})$. This proves the proposition.

LEMMA 6.1. *For the fractional linear vector function (64) there holds*

$$\underline{d}^2 g = \frac{-2\underline{\gamma}^T\underline{h}}{\underline{\gamma}^T\underline{x}+\delta}\underline{d}g \quad and \quad \underline{d}^3 g = 6\left(\frac{\underline{\gamma}^T\underline{h}}{\underline{\gamma}^T\underline{x}+\delta}\right)^2\underline{d}g \qquad (68)$$

*Proof.* We have $(\underline{\gamma}^T\underline{x}+\delta)\underline{g}(\underline{x}) = \underline{\underline{A}}\underline{x}+\underline{\beta}$, whence by the rules in Proposition 3.5

$$(\underline{\gamma}^T\underline{x}+\delta)\underline{\underline{D}}g + \underline{\gamma}^T\otimes\underline{g} = \underline{\underline{A}}$$

and

$$(\underline{\gamma}^T\underline{x}+\delta)\underline{\underline{D}}^2\underline{g} + \underline{\gamma}^T\otimes\underline{\underline{D}}g + \underline{\gamma}^T * \underline{\underline{D}}g = \underline{0} \qquad (69)$$

Via the product rule and the first rule in Proposition 3.4 we get from (69)

$$(\underline{\gamma}^T\underline{x}+\delta)\underline{d}^2\underline{g} + 2(\underline{\gamma}^T\underline{h})\underline{d}g = \underline{0} \qquad (70)$$

and this gives the first expression in (68). Next we use the $\underline{\underline{D}}$ operation on both sides of (70). Using the rules in Proposition 3.5 we first get

$\underline{\underline{D}}(\underline{d}g) = \underline{\underline{D}}^2\underline{g}(\underline{\underline{E}}\otimes\underline{h})$ and $\underline{\underline{D}}(\underline{d}^2\underline{g}) = \underline{\underline{D}}(\underline{\underline{D}}^2\underline{g}(\underline{h}^2)) = \underline{\underline{D}}^3\underline{g}(\underline{\underline{E}}\otimes\underline{h}^2)$ and then from (70)

$$(\underline{\gamma}^T\underline{x}+\delta)\underline{\underline{D}}^3\underline{g}\left(\underline{\underline{E}}\otimes\underline{h}^2\right) + \underline{\gamma}^T\otimes\left(\underline{\underline{D}}^2\underline{g}\,\underline{h}^2\right) + 2(\underline{\gamma}^T\underline{h})\underline{\underline{D}}^2\underline{g}\left(\underline{\underline{E}}\otimes\underline{h}\right) = \underline{0} \qquad (71)$$

Multiplying with $\underline{h}$ from the right on both sides in (71) and using the last rule in Proposition 3.1 with $\underline{\underline{C}} = \underline{\underline{E}}$ we get

$$(\gamma^T \underline{x} + \delta)\underline{\underline{D}}^3 \underline{g} \underline{h}^3 + (\gamma^T \underline{h})(\underline{\underline{D}}^2 \underline{g} \underline{h}^2) + 2(\gamma^T \underline{h})(\underline{\underline{D}}^2 \underline{g} \underline{h}^2) = \underline{0}$$

which is the same as $(\gamma^T \underline{x} + \delta)\underline{d}^3 \underline{g} + 3(\gamma^T \underline{h})\underline{d}^2 \underline{g} = \underline{0}$. Inserting the first expression in (68) herein we get the second expresion. This proves the lemma.

PROPOSITION 6.2. *If the general solution of system (39) has the form (67) then (39) is a vector Riccati equation.*

*Proof.* From Lemma 6.1 follows that

$$\underline{d}^3 \underline{\varphi} \otimes \underline{d}\underline{\varphi} = \underline{d}\underline{\varphi} \otimes \underline{d}^3 \underline{\varphi} = 6\left(\frac{\gamma^T \underline{h}}{\gamma^T \underline{x} + \delta}\right)^2 (\underline{d}\underline{\varphi})^2 \quad , \quad 3(\underline{d}^2 \underline{\varphi})^2 = 12\left(\frac{\gamma^T \underline{h}}{\gamma^T \underline{x} + \delta}\right)^2 (\underline{d}\underline{\varphi})^2$$

The result then follows from Theorem 5.1.

Finally combining Propositions 6.1 and 6.2 we obtain

THEOREM 6.1. *Systen (39) is a vector Riccati equation if and only if its general solution has the form (67).*

As mentioned in the introduction this generalizes a classical result.

## 7. Concluding remarks

REMARK 1. The $*$ notation (30) is introduced only in order to give formula (34) a practical form. However, it is easily seen that there exists one and only one $sn \times sn$ matrix $\underline{\underline{F}}$ which satisfies the condition

$$\forall \underline{u} \in R^s \forall \underline{v} \in R^n : \underline{\underline{F}}(\underline{u} \otimes \underline{v}) = \underline{v} \otimes \underline{u}$$

Then formula (34) can be written

$$\underline{\underline{D}}(\underline{\underline{A}} \otimes \underline{\underline{B}}) = (\underline{\underline{A}} \otimes \underline{\underline{D}}\underline{\underline{B}})(\underline{\underline{F}} \otimes \underline{\underline{E}}_q) + (\underline{\underline{D}}\underline{\underline{A}}) \otimes \underline{\underline{B}}$$

REMARK 2. For the fractional linear vector function (64) there holds

$$\underline{d}\underline{g} \otimes \underline{d}^3 \underline{g} - \frac{3}{2}(\underline{d}^2 \underline{g})^2 = \underline{0} \quad \text{identically}$$

by Lemma 6.1. This can also be written

$$\left[\underline{\underline{D}}\underline{g}(\underline{x}) \otimes \underline{\underline{D}}^3 \underline{g}(\underline{x}) - \frac{3}{2}\left(\underline{\underline{D}}^2 \underline{g}(\underline{x})\right)^{\otimes 2}\right]\underline{h}^4 = \underline{0}$$

identically, i.e., for all $\underline{x} \in R^n$ with $\underline{\gamma}^T \underline{x} + \delta \neq 0$ and for all $\underline{h} \in R^n$. With $\underline{h} = \underline{x}$ this becomes

$$\left[ \underline{D}\underline{g}(\underline{x}) \otimes \underline{D}^3 \underline{g}(\underline{x}) - \frac{3}{2}\left(\underline{D}^2 \underline{g}(\underline{x})\right)^2 \right] \underline{x}^4 = \underline{0} \qquad (72)$$

for all $\underline{x} \in R^n$ with $\underline{\gamma}^T \underline{x} + \delta \neq 0$. It can be shown that if $\Omega \subseteq R^n$ is an open convex set containing $\underline{0}$ and if $\underline{g}(\underline{x}), \underline{x} \in \Omega$ is a $C^3$ function then $\underline{g}(\underline{x})$ has the form (64) if and only if (72) holds for all $\underline{x} \in \Omega$. If n=1 this becomes a classical result.

## 8. Appendix

We prove the propositions stated in section 3.

*Proof of Proposition 3.1.* We have

$$\left(\underline{a}^T \otimes \underline{b}\right)\underline{c} = \left[a_1\underline{b} \quad a_2\underline{b} \quad ... \quad a_n\underline{b}\right]\underline{c} = \sum_{k=1}^{n} a_k \underline{b} c_k = \left(\sum_{k=1}^{n} a_k c_k\right)\underline{b} = \left(\underline{a}^T \underline{c}\right)\underline{b}$$

The r'th block in $\underline{a} \otimes \underline{C}$ is $a_r \underline{C}$, and $\left(a_r \underline{C}\right)\underline{b} = a_r\left(\underline{C}\underline{b}\right)$, whence by the product rule

$$\left(\underline{a} \otimes \underline{C}\right)\underline{b} = \left(\underline{a} \otimes \underline{C}\right)\left((1) \otimes \underline{b}\right) = \underline{a} \otimes \left(\underline{C}\underline{b}\right)$$

In the same way we get

$$\left(\underline{C} \otimes \underline{a}\right)\underline{b} = \left(\underline{C} \otimes \underline{a}\right)\left(\underline{b} \otimes (1)\right) = \left(\underline{C}\underline{b}\right) \otimes \underline{a}$$

This proves the proposition.

*Proof of Proposition 3.2.* It suffices to prove (23) for m=1, i.e., if $\underline{M}$ is a c-symmetric $1 \times n^2$ row $\underline{c}^T = \left[\underline{c}_1^T \quad \underline{c}_2^T \quad ... \quad \underline{c}_n^T\right]$. Let the i'th element in the $1 \times n$ row $\underline{c}_j^T$ be denoted $c_j^i$. We the have $c_j^i = c_i^j$ for all $i, j = 1, 2, ..., n$. Let first p=1, i.e., let $\underline{A}$ be an $n \times 1$ column $\underline{a}$. Then

$$\underline{c}^T\left(\underline{a} \otimes \underline{b}\right) = \left[\underline{c}_1^T \quad \underline{c}_2^T \quad ... \quad \underline{c}_n^T\right]\left[a_1\underline{b} \quad a_2\underline{b} \quad ... \quad a_n\underline{b}\right]^T = \sum_{j=1}^{n} \underline{c}_j^T a_j \underline{b} = \sum_{j,k=1}^{n} c_j^k a_j b_k$$

Interchanging $\underline{a}$ and $\underline{b}$ in this we get

$$\underline{c}^T\left(\underline{b} \otimes \underline{a}\right) = \sum_{j,k=1}^{n} c_j^k b_j a_k = \sum_{j,k=1}^{n} c_j^k a_j b_k = \underline{c}^T\left(\underline{a} \otimes \underline{b}\right) \qquad (73)$$

Next, let $\underline{\underline{A}}$ be an $n \times p$ matrix with columns $\underline{a}_1, \underline{a}_2, ..., \underline{a}_p$. Then

$$\underline{c}^T\left(\underline{\underline{A}} \otimes \underline{b}\right) = \underline{c}^T\left[\underline{a}_1 \otimes \underline{b} \quad \underline{a}_2 \otimes \underline{b} \quad ... \quad \underline{a}_p \otimes \underline{b}\right] = \left[\underline{c}^T\left(\underline{a}_1 \otimes \underline{b}\right) \quad \underline{c}^T\left(\underline{a}_2 \otimes \underline{b}\right) \quad ... \quad \underline{c}^T\left(\underline{a}_p \otimes \underline{b}\right)\right]$$

$$\underline{c}^T\left(\underline{b} \otimes \underline{\underline{A}}\right) = \underline{c}^T\left[\underline{b} \otimes \underline{a}_1 \quad \underline{b} \otimes \underline{a}_2 \quad ... \quad \underline{b} \otimes \underline{a}_p\right] = \left[\underline{c}^T\left(\underline{b} \otimes \underline{a}_1\right) \quad \underline{c}^T\left(\underline{b} \otimes \underline{a}_2\right) \quad ... \quad \underline{c}^T\left(\underline{b} \otimes \underline{a}_p\right)\right]$$

whence $\underline{c}^T\left(\underline{\underline{A}} \otimes \underline{b}\right) = \underline{c}^T\left(\underline{b} \otimes \underline{\underline{A}}\right)$ by (73). This proves the proposition.

*Proof of Proposition 3.3.* If $\underline{a} \otimes \underline{b} + \underline{b} \otimes \underline{a} = \underline{0}$ then $a_j\underline{b} + b_j\underline{a} = \underline{0}$ for all

$j=1,2,...,n$. If $\underline{a} \neq \underline{0}$ then, e.g., $a_i \neq 0$, whence $\underline{b} = -\dfrac{b_i}{a_i}\underline{a}$. Consequently, $\underline{a}$ and $\underline{b}$ Kronecker commute, which implies that $\underline{a} \otimes \underline{b} = \underline{0}$. We infer that

$a_i\underline{b} = \underline{0}$, whence $\underline{b} = \underline{0}$. This proves (27).

If $\underline{a} \otimes \underline{b} + \underline{b} \otimes \underline{a} = \underline{c}^2$ then $a_j\underline{b} + b_j\underline{a} = c_j\underline{c}$ for all $j=1,2,...,n$. Hence also

$a_i\underline{b} + b_i\underline{a} = c_i\underline{c}$ for all $i=1,2,...,n$. This gives $(a_ib_j - a_jb_i)\underline{a} = (a_ic_j - a_jc_i)\underline{c}$

for all $i,j=1,2,...,n$. If there exist $i,j$ such that $a_ic_j - a_jc_i \neq 0$ then

$\underline{c} = \dfrac{a_ib_j - a_jb_i}{a_ic_j - a_jc_i}\underline{a} = \alpha\underline{a}$. Otherwise $a_ic_j - a_jc_i = 0$ for all $i,j=1,2,...,n$, whence

$a_i\underline{c} = c_i\underline{a}$ for all $i=1,2,...,n$. If $\underline{a} \neq \underline{0}$ then, e.g., $a_k \neq 0$, giving $\underline{c} = \dfrac{c_k}{a_k}\underline{a} = \alpha\underline{a}$.

Otherwise $\underline{a} = \underline{0}$, whence $\underline{c}^2 = \underline{0}$, i.e., $\underline{c} = \underline{0} = 1 \cdot \underline{a}$. In this argumentation we may interchange $\underline{a}$ and $\underline{b}$, giving $\underline{c} = \beta\underline{b}$. This proves (28).

Suppose now that $\underline{a} \neq \underline{0}$. If $\beta \neq 0$ in (28), then $\underline{b} = \dfrac{\alpha}{\beta}\underline{a} = \gamma\underline{a}$. If $\beta = 0$ then

$\underline{c} = \underline{0}$, whence $\underline{b} = \underline{0}$ by (27), which implies that $\underline{b} = 0 \cdot \underline{a} = \gamma\underline{a}$. Inserting

$\underline{b} = \gamma\underline{a}$ and $\underline{c} = \alpha\underline{a}$ we get $2\gamma\underline{a} \otimes \underline{a} = \alpha^2\underline{a} \otimes \underline{a}$. Since $\underline{a} \neq \underline{0}$ we infer that

$2\gamma = \alpha^2$, whence $\underline{b} = \dfrac{1}{2}\alpha^2\underline{a}$. This proves (29).

LEMMA 7.1. *Let $\underline{\underline{A}}, \underline{\underline{A}}_1, ..., \underline{\underline{A}}_s$ be $m \times n$ matrices and let $\underline{\underline{B}}, \underline{\underline{B}}_1, ..., \underline{\underline{B}}_s$ be $n \times p$*

*matrices. Then*

$$\underline{\underline{A}}\left[\underline{\underline{B}}_1 \quad \underline{\underline{B}}_2 \quad ... \quad \underline{\underline{B}}_s\right] = \left[\underline{\underline{A}}\underline{\underline{B}}_1 \quad \underline{\underline{A}}\underline{\underline{B}}_2 \quad ... \quad \underline{\underline{A}}\underline{\underline{B}}_s\right] \tag{74}$$

and

$$\left[\underline{\underline{A}}_1 \quad \underline{\underline{A}}_2 \quad ... \quad \underline{\underline{A}}_s\right] \otimes \underline{\underline{B}} = \left[\underline{\underline{A}}_1 \otimes \underline{\underline{B}} \quad \underline{\underline{A}}_2 \otimes \underline{\underline{B}} \quad ... \quad \underline{\underline{A}}_s \otimes \underline{\underline{B}}\right] \tag{75}$$

and

$$\left[\underline{\underline{A}}_1 \quad \underline{\underline{A}}_2 \quad ... \quad \underline{\underline{A}}_s\right]\left(\underline{\underline{E}}_s \otimes \underline{\underline{B}}\right) = \left[\underline{\underline{A}}_1\underline{\underline{B}} \quad \underline{\underline{A}}_2\underline{\underline{B}} \quad ... \quad \underline{\underline{A}}_s\underline{\underline{B}}\right] \tag{76}$$

*Proof.* The rules (74) and (75) are obvious. Let $\underline{a}_{i,j}^T$ denote the j'th row in $\underline{A}_i$, $i=1,2,...,s$ , $j=1,2,...,m$. Then the j'th row of the two sides in (76) are both $\begin{bmatrix} \underline{a}_{1,j}^T\underline{B} & \underline{a}_{2,j}^T\underline{B} & ... & \underline{a}_{s,j}^T\underline{B} \end{bmatrix}$. This proves (76).

LEMMA 7.2. *Let* $\underline{A}_1,\underline{A}_2,...,\underline{A}_p$ *be* $m \times n$ *matrices and let* $\underline{B}$ *be a* $p \times q$ *matrix. Then*

$$\begin{bmatrix} \underline{A}_1 & \underline{A}_2 & ... & \underline{A}_p \end{bmatrix}(\underline{B} \otimes \underline{E}_n) = \sum_{k=1}^{p}\begin{bmatrix} b_{k1}\underline{A}_k & b_{k2}\underline{A}_k & ... & b_{kq}\underline{A}_k \end{bmatrix} \qquad (77)$$

*Proof.* Let $\underline{a}_{i,j}^T$ denote the j'th row in $\underline{A}_i$, $i=1,2,...,p$ , $j=1,2,...,m$. These are all $1 \times n$ rows. Let $\underline{b}_i$, $i=1,2,...,q$ denote the i'th column in $\underline{B}$ , all $p \times 1$ columns. The product of the j'th row in $\begin{bmatrix} \underline{A}_1 & \underline{A}_2 & ... & \underline{A}_p \end{bmatrix}$ and the i'th block of columns in $\underline{B} \times \underline{E}_n$ then is

$$\begin{bmatrix} \underline{a}_{1,j}^T & \underline{a}_{2,j}^T & ... & \underline{a}_{p,j}^T \end{bmatrix}\begin{bmatrix} b_{1i}\underline{E}_n & b_{2i}\underline{E}_n & ... & b_{pi}\underline{E}_n \end{bmatrix}^T = \sum_{k=1}^{p}\underline{a}_{k,j}^T b_{ki}$$

i=1,2,...,q , j=1,2,...,p. Hence, the j'th row on the left hand side of (77) is

$$\begin{bmatrix} \sum_{k=1}^{p}b_{k1}\underline{a}_{k,j}^T & \sum_{k=1}^{p}b_{k2}\underline{a}_{k,j}^T & ... & \sum_{k=1}^{p}b_{kq}\underline{a}_{k,j}^T \end{bmatrix} = \sum_{k=1}^{p}\begin{bmatrix} b_{k1}\underline{a}_{k,j}^T & b_{k2}\underline{a}_{k,j}^T & ... & b_{kq}\underline{a}_{k,j}^T \end{bmatrix}$$

which is the j'th row in the matrix on the right hand side of (77). This proves the lemma.

*Proof of Proposition 3.4.* Let $\underline{C}_i$ denote the i'th column in $\underline{C}$. Then

$$(\underline{A} * \underline{C})\underline{h}^2 = \begin{bmatrix} \underline{A} \otimes \underline{C}_1 & \underline{A} \otimes \underline{C}_2 & ... & \underline{A} \otimes \underline{C}_n \end{bmatrix}\underline{h}^2 = \sum_{k=1}^{n}h_k(\underline{A} \otimes \underline{C}_k)\underline{h}$$

Using the last rule in (23) i Proposition 3.1 we further get

$$(\underline{A} * \underline{C})\underline{h}^2 = \sum_{k=1}^{n}h_k((\underline{A}\underline{h}) \otimes \underline{C}_k) = (\underline{A}\underline{h}) \otimes \sum_{k=1}^{n}h_k\underline{C}_k = (\underline{A}\underline{h}) \otimes (\underline{C}\underline{h})$$

which is the first rule in (32).By the product rule we get

$$(\underline{A} * \underline{B})\underline{h}^3 = \begin{bmatrix} \underline{A} \otimes \underline{B}_1 & \underline{A} \otimes \underline{B}_2 & ... & \underline{A} \otimes \underline{B}_n \end{bmatrix}\begin{bmatrix} h_1\underline{h}^2 & h_2\underline{h}^2 & ... & h_n\underline{h}^2 \end{bmatrix}^T$$

$$= \sum_{k=1}^{n}h_k(\underline{A} \otimes \underline{B}_k)\underline{h}^2 = \sum_{k=1}^{n}h_k((\underline{A}\underline{h}) \otimes (\underline{B}_k\underline{h}))$$

$$= (\underline{A}\underline{h}) \otimes \sum_{k=1}^{n}h_k(\underline{B}_k\underline{h}) = (\underline{A}\underline{h}) \otimes (\underline{B}\underline{h}^2)$$

which is the second rule.

*Proof of Proposition 3.5.* We have

$$\frac{\partial}{\partial x_j}\left(\underline{A}\,\underline{B}\right)=\underline{A}\,\frac{\partial}{\partial x_j}\,\underline{B}+\left(\frac{\partial}{\partial x_j}\,\underline{A}\right)\underline{B}$$

j=1,2,...,s, which may be collected to (33) by means of the rules (74) and (76) in Lemma 7.1. Furthermore,

$$\frac{\partial}{\partial x_j}\left(\underline{A}\otimes\underline{B}\right)=\underline{A}\otimes\frac{\partial}{\partial x_j}\,\underline{B}+\left(\frac{\partial}{\partial x_j}\,\underline{A}\right)\otimes\underline{B}$$

j=1,2,...,s, which may be collected to (34) by means of the rule (74) in Lemma 7.1 and the definition (30). This proves the proposition.

*Proof of Proposition 3.6.* Let the variables in $R^s$ be denoted $\underline{\xi}$, and let the j'th column in $\underline{A}(\underline{x})$ be denoted $\underline{a}_j(\underline{x})$ (an $m\times 1$ column). By means of the differentiation rule for composite functions we get omitting all arguments

$$\frac{\partial}{\partial\xi_j}\left(\underline{A}\circ\underline{\varphi}\right)=\frac{\partial}{\partial\xi_j}\left[\underline{a}_1\circ\underline{\varphi}\quad\underline{a}_2\circ\underline{\varphi}\quad\ldots\quad\underline{a}_n\circ\underline{\varphi}\right]$$

$$=\left[\sum_{k=1}^{p}\left(\frac{\partial\underline{a}_1}{\partial x_k}\circ\underline{\varphi}\right)\frac{\partial\varphi_k}{\partial\xi_j}\quad\sum_{k=1}^{p}\left(\frac{\partial\underline{a}_2}{\partial x_k}\circ\underline{\varphi}\right)\frac{\partial\varphi_k}{\partial\xi_j}\quad\ldots\quad\sum_{k=1}^{p}\left(\frac{\partial\underline{a}_n}{\partial x_k}\circ\underline{\varphi}\right)\frac{\partial\varphi_k}{\partial\xi_j}\right]$$

$$=\sum_{k=1}^{p}\left(\frac{\partial\underline{A}}{\partial x_k}\circ\underline{\varphi}\right)\frac{\partial\varphi_k}{\partial\xi_j}$$

j=1,2,...,s. Via the rule in Lemma 7.2 this can be collected to

$$\underline{D}\left(\underline{A}\circ\underline{\varphi}\right)=\left[\sum_{k=1}^{p}\left(\frac{\partial\underline{A}}{\partial x_k}\circ\underline{\varphi}\right)\frac{\partial\varphi_k}{\partial\xi_1}\quad\sum_{k=1}^{p}\left(\frac{\partial\underline{A}}{\partial x_k}\circ\underline{\varphi}\right)\frac{\partial\varphi_k}{\partial\xi_2}\quad\ldots\quad\sum_{k=1}^{p}\left(\frac{\partial\underline{A}}{\partial x_k}\circ\underline{\varphi}\right)\frac{\partial\varphi_k}{\partial\xi_s}\right]$$

$$=\sum_{k=1}^{p}\left[\frac{\partial\varphi_k}{\partial\xi_1}\left(\frac{\partial\underline{A}}{\partial x_k}\circ\underline{\varphi}\right)\quad\frac{\partial\varphi_k}{\partial\xi_2}\left(\frac{\partial\underline{A}}{\partial x_k}\circ\underline{\varphi}\right)\quad\ldots\quad\frac{\partial\varphi_k}{\partial\xi_s}\left(\frac{\partial\underline{A}}{\partial x_k}\circ\underline{\varphi}\right)\right]$$

$$=\left[\left(\underline{D}\underline{A}\right)\circ\underline{\varphi}\right]\left[\left(\underline{D}\underline{\varphi}\right)\otimes\underline{E}_n\right]$$

This proves the proposition.

*Proof of Proposition 3.7.* Let $\underline{c}_{i,j}$ denote the j'th column in $\underline{\underline{C}}_i$, i,j=1,2,...,n. We put $\alpha_{i,j}=\underline{e}_n^T\underline{c}_{i,j}$ (this is the n'th element in $\underline{c}_{i,j}$). Then

$$\alpha(\underline{h})h_n=\underline{e}_n^T\underline{\underline{C}}\underline{h}^2=\underline{e}_n^T\left[\underline{\underline{C}}_1\quad\underline{\underline{C}}_2\quad\ldots\quad\underline{\underline{C}}_n\right]\left[h_1\underline{h}\quad h_2\underline{h}\quad\ldots\quad h_n\underline{h}\right]^T$$

$$=\underline{e}_n^T\sum_{i=1}^{n}h_i\underline{\underline{C}}_i\underline{h}=\underline{e}_n^T\sum_{i,j=1}^{n}h_ih_j\underline{c}_{i,j}=\sum_{i,j=1}^{n}\alpha_{i,j}h_ih_j=\tag{78}$$

$$=\alpha_{n,n}h_n^2+\left(\sum_{i=1}^{n-1}(\alpha_{n,i}+\alpha_{n,i})h_i\right)h_n+\sum_{i,j=1}^{n-1}\alpha_{i,j}h_ih_j$$

for all $\underline{h} \in R^n$. Putting $h_n = 0$ we obtain that the last sum on the right hand side of (78) is zero for all $h_1, h_2, ..., h_{n-1}$. Then we infer that

$$\alpha(\underline{h}) = \alpha_{n,n} h_n + \sum_{i=1}^{n-1} (\alpha_{n,i} + \alpha_{i,n}) h_i$$

for $h_n \neq 0$. This can be written $\alpha(\underline{h}) = \underline{a}^T \underline{h}$. Since $\alpha(\underline{h}) = \dfrac{1}{h_i} \underline{e}_i^T \underline{\underline{C}} \underline{h}^2$ if $h_i \neq 0$ then $\alpha(\underline{h})$ is a continuous function for $\underline{h} \neq \underline{0}$. We infer that $\alpha(\underline{h}) = \underline{a}^T \underline{h}$ for all $\underline{h} \in R^n$. From the product rule we get

$$\left( \underline{a}^T \otimes \underline{\underline{E}} + \underline{\underline{E}} \otimes \underline{a}^T \right) \underline{h}^2 = (\underline{a}^T \underline{h}) \otimes \underline{h} + \underline{h} \otimes (\underline{a}^T \underline{h}) = 2(\underline{a}^T \underline{h}) \underline{h}$$

whence $\underline{\underline{C}} \underline{h}^2 = \dfrac{1}{2} \left( \underline{a}^T \otimes \underline{\underline{E}} + \underline{\underline{E}} \otimes \underline{a}^T \right) \underline{h}^2$. By Example 2.1 the matrix $\dfrac{1}{2} \left( \underline{a}^T \otimes \underline{\underline{E}} + \underline{\underline{E}} \otimes \underline{a}^T \right)$ is c-symmetric. So, if also $\underline{\underline{C}}$ is also c-symmetric we infer that $\underline{\underline{C}} = \dfrac{1}{2} \left( \underline{a}^T \otimes \underline{\underline{E}} + \underline{\underline{E}} \otimes \underline{a}^T \right)$ by Corollary 2.2. This proves the proposition.

Department of Mathematics, Technical University of Denmark,

e-mail: K.M.Andersen@mat.dtu.dk